# Counting the Number of Non-Equivalent Classes of Fuzzy Matrices Using Combinatorial Techniques

S. R. Kannan [1] and Rajesh Kumar Mohapatra [1, *]
[1]Pondicherry University (A Central University of India), Puducherry, India

*Corresponding author: mohapatrarajesh030@gmail.com

**Abstract.** The novelty of this paper is to construct the explicit combinatorial formula for the number of all distinct fuzzy matrices of finite order, which leads us to invent a new sequence. In order to achieve this new sequence, we analyze the behavioral study of equivalence classes on the set of all fuzzy matrices of a given order under a suitable natural equivalence relation. In addition this paper characterizes the properties of non-equivalent classes of fuzzy matrices of order $n$ with elements having degrees of membership values anywhere in the closed unit interval [0,1]. Further, this paper also derives some important relevant results by enumerating the number of all distinct fuzzy matrices of a given order in general. And also, we achieve these results by incorporating the notion of $k$-level fuzzy matrices, chains, and flags (maximal chains).

## 1.Introduction

One of the most challenging problems of fuzzy matrix theory is to classify for concerning the study of the set of all fuzzy matrices of any given order and to count the same. Counting of all distinct fuzzy matrices of finite order is great of an interest in both and theoretical and practical point of view in contexts of mathematics. In recent years, this topic has undergone a remarkable development in many areas and many interesting important results have been proposed. The pioneering work of Zadeh [1] on fuzzy subsets of a set and Rosenfeld [2] on fuzzy subgroups of a group led to the fuzzification of some algebraic structures like group, ring, etc. Though handling such a problem is quite a difficult job, many researchers have been still working on classifying and counting both fuzzy subsets and fuzzy subgroups in the last few years. Research works on counting the number of fuzzy subsets of a finite set was initially noticed by Murali and Makamba [3] and later many other researchers such as Tărnăucean [4], Šešelja and Tepavčević [5], and Jain [6] have taken further away. To determine the number of fuzzy subgroups of finite groups, at first, several papers have treated the special cases of finite abelian groups such as computing the number of distinct square-free order of fuzzy cyclic group [3], fuzzy cyclic group of order $p^n q^m$ ($p, q$ primes) ([7], [8]–[11]). Later, the authors in ([4], [12]) deal for determining the number of distinct fuzzy subgroups for two classes of finite abelian groups: finite cyclic groups and finite elementary abelian $p$-groups. Further, this investigation has been propagated to some remarkable classes of non-abelian groups: symmetric groups ([13], [14]), alternating group [15], dihedral groups ([16], [17], [18]), hamiltonian groups [19]. Subsequently, the same problems were also implemented and analyzed in the case of the fuzzy normal subgroup ([20], [16]).

The research works on calculating the formula for the number of all distinct fuzzy matrices through a combinatorial approach are very important in deriving a new sequence that can be used in many real-life applications. Researchers have concentrated their research activities on counting the fuzzy subsets of a finite set as well as fuzzy subgroups of a finite group, but no researchers have focused their research works on computing the number of fuzzy matrices. Since the problem of counting the number of non-equivalent fuzzy subsets in one-dimensional space was obtained, but it remains still open for two-dimensional space, i.e., the formula for calculating the number of non-equivalent classes of fuzzy matrices. Thus, it motivates us to contribute novel research for invention in the field of fuzzy matrix theory ([21], [22]). In order to achieve the new invents on establishing the formula for the number of all distinct fuzzy matrices on finite set, this paper discusses equivalence classes on fuzzy matrices on a finite set under a natural equivalence relation. Because of the fact that fuzzy matrices on a finite set are more abundant than crisp matrices. The necessity of classification of fuzzy matrices have done based on $\alpha$-cuts due to the following fact: there are uncountably many fuzzy matrices on the same domain, even if the domain is either a singleton (or finite) or countable. Since there are different notions of equivalence ([4], [5], [11], [23]–[27]), research works have grown to be a major branch which have many more interesting challenges in the classification of fuzzy subsets, fuzzy subgroups, etc. The present paper utilizes equivalence [24] by accepting the property on support which is significant to construct the formula for the number of all distinct fuzzy matrices. In this paper we mainly focus to propose a closed explicit formula by indicating the number of all distinct fuzzy matrices of order $n$. This closed formula leads us to invent a new important sequence that is not present in the Online Encyclopedia of Integer Sequences (OEIS) [28].

The structure of this paper is organized as follows: Section 2 focuses on a few important definitions, results and notational set up which are necessary throughout the paper. Precisely, we will concentrate on counting the number of all fuzzy matrices up to Murali's equivalence relation and will derive some of its relevant results in Section 3. In the final section conclusions and further research of this paper are specified.

## 2. Preliminaries

The main goal of this section is to set up notations, collect some basic definitions and results with its properties to introduce the new important results on counting some specific fuzzy matrices.

### 2.1. Fuzzy matrices (FMs)

Throughout this paper, let us assume that $\tilde{A} = (a_{ij})$ be an $n \times n$ fuzzy matrix with elements having membership values in the real unit interval $I = [0,1]$, where $n$ is a non-zero, non-negative integer. The *union* ($\cup$), *intersection* ($\cap$) of two fuzzy matrices, and *complementation* ($^c$) of a fuzzy matrix are defined by using supremum (sup or max) and infimum (inf or min) component-wise, and $1 - a_{ij}$ operator pointwise, respectively [1], [29]. The *containment* of a fuzzy matrix $\tilde{A}$ in a fuzzy matrix $\tilde{B}$, denoted as $\tilde{A} \subseteq \tilde{B}$ if $a_{ij} \leq b_{ij}$ for all $i,j$. Further we denote fuzzy matrices by $\tilde{A}, \tilde{B}, \tilde{C}$, etc. and its membership values by $\alpha, \beta, \gamma$, etc. Through an $\alpha$-*cut* of fuzzy matrix $\tilde{A}$ for $\alpha$ belongs to $I$, we have a crisp matrix $\tilde{A}^\alpha = \begin{cases} 1, & a_{ij} \geq \alpha \\ 0, & \text{otherwise} \end{cases}$. This is called the *weak $\alpha$-cut*. By *strong $\alpha$-cut* we mean $\tilde{A}_\alpha = \{1, \text{if } a_{ij} > \alpha \text{ and } 0, \text{otherwise}\}$. But, in this paper we are always dealt with weak $\alpha$-cut. It is easy to verify that for $0 \leq \alpha \leq \alpha' \leq 1$, we have $\tilde{A}^\alpha \supseteq \tilde{A}^{\alpha'}$. For any arbitrary fuzzy matrix $\tilde{A}$ it can be decomposed into a union of characteristic function as follows:

**Theorem 2.1** [30], [29], [31]. *For any fuzzy matrix $\tilde{A} = \underset{\alpha}{\cup}\{\alpha \chi_{\tilde{A}^\alpha} : 0 \leq \alpha \leq 1\}$, where $\chi_{\tilde{A}^\alpha}$ denotes the characteristic function of the crisp matrix $\tilde{A}^\alpha$.*

## 2.2. Binomial numbers

To prove our results, we recall some known definitions and theorems in the area of combinatorics.

**Definition 2.2** (*Binomial Number*) [32]**.** Let us denote $\binom{n}{m}$ be the number of $m$-element subsets of an $n$-element set; that is nothing but the number of ways we can select $m$ distinct elements from an $n$-element set. This is well-known as a *binomial number* or a *binomial coefficient*. (*An alternative notation,* $C(n,m)$)

Now we are going to present a famous theorem, known as binomial theorem.
**Theorem 2.3** (*Binomial Theorem*) [32]**.** For all integers $n \geq 0$,

$$(a+b)^n = \sum_{m=0}^{n} \binom{n}{m} a^m b^{n-m}.$$

In fact, one has to give the familiar formula of binomial numbers. (By definition, $0! = 1$.)
**Lemma 2.4.** For $n \geq m \geq 0$,

$$\binom{n}{m} = \frac{n!}{(n-m)!},$$

with the convention that $\binom{n}{m} = 0$, for any $m > n$.

## 2.3. Equivalent fuzzy matrices and concepts of chains

The purpose of this section is to briefly discuss the study of an equivalence relation on the set of all fuzzy matrices and the concept of chains. Thus we start with an equivalence relation $\cong$ defined on any class of fuzzy matrices as follows:

**Definition 2.5** [24], [33], [34]**.** $\tilde{A} \cong \tilde{B}$ if and only if for all $i, j, r, s$

(i) $a_{ij} > a_{rs}$ if and only if $b_{ij} > b_{rs}$
(ii) $a_{ij} = 1$ if and only if $b_{ij} = 1$
(iii) $a_{ij} = 0$ if and only if $b_{ij} = 0$.

It is easy to verify that this relation is indeed an equivalence relation on the set of all fuzzy matrices and when its entries restricted to $I' = \{0, 1\}$, which corresponds with equality of crisp matrices. Based on this equivalence relation, equivalence class having $\tilde{A}$ is denoted as $[\tilde{A}]$ and two fuzzy matrices $\tilde{A}$ and $\tilde{B}$ are distinct if $\tilde{A}$ and $\tilde{B}$ are not equivalent, that is, $\tilde{A} \not\cong \tilde{B}$.

The next proposition proposes the relation between $\alpha$-cuts and equivalence.
**Proposition 2.6** [35]**.** *Let $\tilde{A}$ and $\tilde{B}$ be two fuzzy matrices of order n. Then $\tilde{A} \cong \tilde{B}$ if and only if for each $\alpha > 0$ there exists a $\beta > 0$ such that $\tilde{A}^\alpha = \tilde{B}^\beta$.*

It follows from the above proposition that the equivalent fuzzy matrices can be analyzed by their $\alpha$-cuts. This observation promotes us an innovation to raise the ideas $k$-level fuzzy matrices, chains, and flags. Now we will elaborate them explicitly in the following subsection.

## 2.4. $k$-level fuzzy matrices and flags

We will classify non-equivalent classes of fuzzy matrices by incorporating the concept $\alpha$-cuts for $0 \leq \alpha \leq 1$. In this subsection, we briefly give some definitions using the basic ideas.

**Definition 2.7.** For $n \in \mathbb{N}$ and $0 \leq k \leq n^2$, by a $k$-level fuzzy matrix $\tilde{A}$, it means $\tilde{A}$ has $k$-number of distinct membership values in the open unit interval, that is, $I \setminus \{0,1\}$. Explicitly, a $k$-level fuzzy matrix of order $n$ is a $(k+1)$-pair of a chain of crisp matrices under the usual inclusion of the form

$$\Lambda : \tilde{A}^{\alpha_0} \subset \tilde{A}^{\alpha_1} \subset \tilde{A}^{\alpha_2} \subset \cdots \subset \tilde{A}^{\alpha_k}$$

with

$$1 \geq \alpha_0 > \alpha_1 > \alpha_2 > \cdots > \alpha_k \geq 0$$

$\alpha_i$'s are in the unit interval $I$, and not necessarily including 0 and 1, written in the descending order of magnitude. Then a fuzzy matrix $\tilde{A} = \bigcup_{\alpha_i} \{\alpha_i \chi_{\tilde{A}^{\alpha_i}} : 0 \leq \alpha_i \leq 1\}$ has its $\alpha_i$-cuts equals $\tilde{A}^{\alpha_i}$.

Here $\tilde{A}^{\alpha_i}$'s are called various components of the chain $\Lambda$. It follows that the above inclusions in the chain $\Lambda$ are always taken to be strict. By taking $k = n^2$ in the above chain $\Lambda$, we can get a maximal chain, that is named as a *flag*.

We call two $k$-level fuzzy matrices are distinct if they are not equivalent. It is also clear that for any fuzzy matrix of order $n$ can have maximum $n^2$-distinct values of membership degree in $I$. Hence the maximum possible number of distinct $\alpha$-cut relational matrices of a fuzzy matrix of order $n$ is $n^2 + 1$. The cause behind is that the set of $n^2$-distinct real numbers in the open interval $(0,1)$ will split the closed interval $[0,1]$ into $n^2 + 1$ segments. Now take an $\alpha$ in each open segment will provide a distinct $\alpha$-cut. Therefore, the length of the flag for a fuzzy matrix of order $n$ is $n^2 + 1$.

Hence, it can be concluded that $\tilde{A} \cong \tilde{B}$ if and only if $\tilde{A}$ and $\tilde{B}$ determine the same chain of crisp matrices of type $\Lambda$. It obviously follows that there is a one-to-one correspondence between the set of $k$-level distinct fuzzy matrices of order $n$ and the set of chains crisp matrices of order $n$ of length $k$ under the usual inclusion.

**Definition 2.8.** Let $n$ be a positive integer and let $A_0, A_1, A_2, \ldots A_k$ be crisp matrices of order $n$ with $A_i \subset A_{i+1}$ for $i = \overline{0, k-1}$ with $0 \leq k \leq n^2$ then the following matrix chain of the type

$$\Gamma : A_0 \subset A_1 \subset A_2 \subset \cdots \subset A_k$$

is called a chain of crisp matrices of order $n$. In this case, the integer $k$ is called the length of the proper subgroup chain $\Gamma$, and the subgroups $A_0$ and $A_k$ are called the initial term and the terminal term of $\Gamma$.

**Definition 2.9.** (i) Let $O$ be an $n \times n$ null matrix is a matrix that is defined as $\tilde{A} = (a_{ij})_{n \times n}$, where $a_{ij} = 0$ for all $i, j \in \{1, 2, \ldots, n\}$ and $J$ be a unit matrix of order $n$ is a matrix which is defined as $J = (a_{ij})_{n \times n}$, where $a_{ij} = 1$ for every $1 \leq i, j \leq n$.

(ii) A chain of crisp matrices of order $n$ is called $O$-rooted (respectively, $J$-rooted) if it contains $O$ (respectively, $J$). Otherwise, it is simply called a chain of crisp matrices of order $n$.

It can be remark that there is a one-to-one correspondence between the set of $O$-rooted (respectively, $J$-rooted) $k$-level distinct fuzzy matrices of order $n$ and the set of $O$-rooted (respectively, $J$-rooted) chains of crisp of matrices of order $n$ of length $k$ under inclusion.

**Definition 2.10.** Suppose $n$ be a positive integer and $0 \leq k \leq n^2$. Let $\mathcal{M}$ be the set of all crisp matrices of order $n$ and let $C$ be the family of all chains of crisp matrices of order $n$ of length $k$. Now set

$FM_{n,k}(\mathcal{M}) = \{\Gamma \in C \mid$ the length of $\Gamma$ is $k$, that is the initial term and the terminal term of $\Gamma$ are not necessarily be the null matrix and the unit matrix respectively$\}$.

$FM_{n,k}^O(\mathcal{M}) = \{\Gamma \in C \mid$ the initial term of $\Gamma$ is $A_0 = O$ of length $k\}$;

$FM_{n,k}^J(\mathcal{M}) = \{\Gamma \in C \mid$ the terminal term of $\Gamma$ is $A_k = J$ of length $k\}$;

Let us denote $FM_{n,k}(\mathcal{M})$, $FM_{n,k}^O(\mathcal{M})$ and $FM_{n,k}^J(\mathcal{M})$ be the set of all chains of crisp matrices of order $n$ of length $k$, $O$-rooted chains of crisp matrices of order $n$ of length $k$ and $J$-rooted chains of crisp matrices of order $n$ of length $k$, respectively. We use notation $f_{n,k}$, $f_{n,k}^O$ and $f_{n,k}^J$ to denote the cardinal numbers of $FM_{n,k}(\mathcal{M})$, $FM_{n,k}^O(\mathcal{M})$ and $FM_{n,k}^J(\mathcal{M})$ respectively. The relations between these numbers are of the following.

**Remark 2.11.** It is obvious that $FM_{n,k}^O(\mathcal{M}) \subset FM_{n,k}(\mathcal{M})$ and $FM_{n,k}^J(\mathcal{M}) \subset FM_{n,k}(\mathcal{M})$. Hence, $f_{n,k}^O < f_{n,k}$ and $f_{n,k}^J < f_{n,k}$.

Also, we denote $FM_n(\mathcal{M})$, $FM_n^O(\mathcal{M})$ and $FM_n^J(\mathcal{M})$ be the set of all chains of crisp matrices of order $n$, $O$-rooted chains of crisp matrices of order $n$ and $J$-rooted chains of crisp matrices of order $n$, respectively. Take $f_n = |FM_n(\mathcal{M})|$, $f_n^O = |FM_n^O(\mathcal{M})|$ and $f_n^J = |FM_n^J(\mathcal{M})|$.

Similarly, we can conclude that there is also a one-to-one correspondence between the collection of distinct fuzzy matrices of order $n$ (respectively, $O$-rooted distinct fuzzy matrices of order $n$, $J$-rooted distinct fuzzy matrices of order $n$) and the collection of chains crisp matrices of order $n$ (respectively, $O$-rooted chains crisp matrices of order $n$, $J$-rooted chains crisp matrices of order $n$) under inclusion.

### 3. Enumeration of fuzzy matrices

In this section, we keen to develop the method for counting the number of $FM(\mathcal{M})$. Here, the explicit closed combinatorial summation formulae were discovered for computing the number of $FM_n(\mathcal{M})$, $FM_n^O(\mathcal{M})$ and $FM_n^J(\mathcal{M})$. In order to determine these numbers, we first establish the number of $FM_{n,k}(\mathcal{M})$, $FM_{n,k}^O(\mathcal{M})$ and $FM_{n,k}^J(\mathcal{M})$ in the following section.

## 3.1. The number of k-level fuzzy matrices

This subsection motivates to achieve our objective of finding the formulas for calculating the numbers $f_{n,k}, f^O_{n,k}$ and $f^J_{n,k}$. In order to derive the explicit formula of $f_{n,k}$, we shall first investigate with some of its specific cases.

Let us take $n = 2$, the lattice of $FM_{2,0}(\mathcal{M})$ (that is the set of all crisp matrices of order 2), $L(FM_{2,0}(\mathcal{M}))$ is constituted by the following matrices:

$$A_0 = \begin{pmatrix} 0 & 0 \\ 0 & 0 \end{pmatrix},$$

$$A^1_1 = \begin{pmatrix} 1 & 0 \\ 0 & 0 \end{pmatrix}, A^2_1 = \begin{pmatrix} 0 & 1 \\ 0 & 0 \end{pmatrix}, A^3_1 = \begin{pmatrix} 0 & 0 \\ 1 & 0 \end{pmatrix}, A^4_1 = \begin{pmatrix} 0 & 0 \\ 0 & 1 \end{pmatrix},$$

$$A^{1,2}_2 = \begin{pmatrix} 1 & 1 \\ 0 & 0 \end{pmatrix}, A^{1,3}_2 = \begin{pmatrix} 1 & 0 \\ 1 & 0 \end{pmatrix}, A^{1,4}_2 = \begin{pmatrix} 1 & 0 \\ 0 & 1 \end{pmatrix}, A^{2,3}_2 = \begin{pmatrix} 0 & 1 \\ 1 & 0 \end{pmatrix}, A^{2,4}_2 = \begin{pmatrix} 0 & 1 \\ 0 & 1 \end{pmatrix},$$
$$A^{3,4}_2 = \begin{pmatrix} 0 & 0 \\ 1 & 1 \end{pmatrix},$$

$$A^{1,2,3}_3 = \begin{pmatrix} 1 & 1 \\ 1 & 0 \end{pmatrix}, A^{1,2,4}_3 = \begin{pmatrix} 1 & 1 \\ 0 & 1 \end{pmatrix}, A^{1,3,4}_3 = \begin{pmatrix} 1 & 0 \\ 1 & 1 \end{pmatrix}, A^{2,3,4}_3 = \begin{pmatrix} 0 & 1 \\ 1 & 1 \end{pmatrix},$$

$$A_4 = \begin{pmatrix} 1 & 1 \\ 1 & 1 \end{pmatrix},$$

and its lattice structure has shown in Figure 1.

For $n = 2$ and $k = 3$. We can precisely determine the number of $FM_{2,3}(\mathcal{M})$, by describing all chains through manually by direct calculation in the five possibilities as listed below:

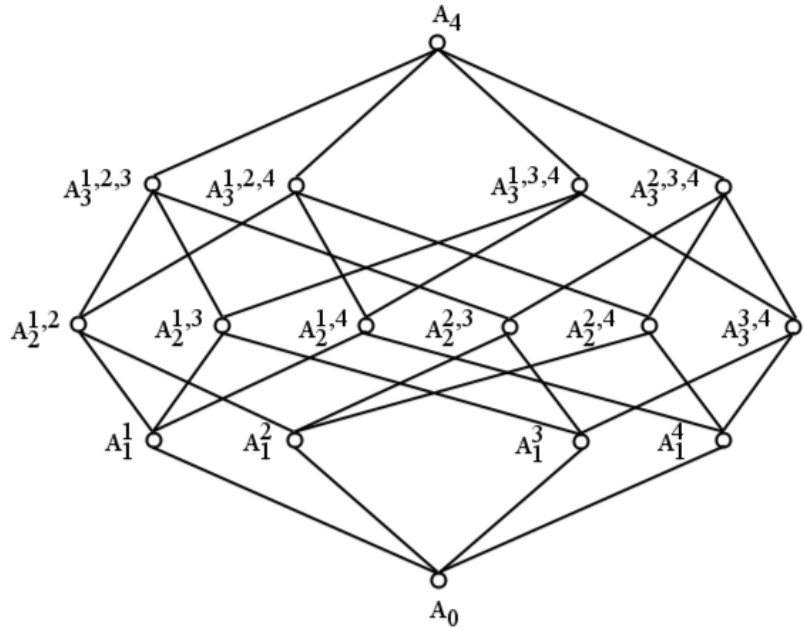

**Figure 1.** Graphical illustration of $L(FM_{2,0}(\mathcal{M}))$

**Case I.** Consider the chains in $FM_{2,3}(\mathcal{M})$ of the type

$$A_0 \subset A_1 \subset A_2 \subset A_3$$

have the following 24 chains:

$A_0 \subset A_1^1 \subset A_2^{1,2} \subset A_3^{1,2,3}, A_0 \subset A_1^1 \subset A_2^{1,2} \subset A_3^{1,2,4}, A_0 \subset A_1^1 \subset A_2^{1,3} \subset A_3^{1,2,3}, A_0 \subset A_1^1$
$\subset A_2^{1,3} \subset A_3^{1,3,4}, A_0 \subset A_1^1 \subset A_2^{1,4} \subset A_3^{1,2,4}, A_0 \subset A_1^1 \subset A_2^{1,4} \subset A_3^{1,3,4},$

$A_0 \subset A_1^2 \subset A_2^{1,2} \subset A_3^{1,2,3}, A_0 \subset A_1^2 \subset A_2^{1,2} \subset A_3^{1,2,4}, A_0 \subset A_1^2 \subset A_2^{2,3} \subset A_3^{1,2,3}, A_0 \subset A_1^2$
$\subset A_2^{2,3} \subset A_3^{2,3,4}, A_0 \subset A_1^2 \subset A_2^{2,4} \subset A_3^{1,2,4}, A_0 \subset A_1^2 \subset A_2^{2,4} \subset A_3^{2,3,4},$

$A_0 \subset A_1^3 \subset A_2^{1,3} \subset A_3^{1,2,3}, A_0 \subset A_1^3 \subset A_2^{1,3} \subset A_3^{1,3,4}, A_0 \subset A_1^3 \subset A_2^{2,3} \subset A_3^{1,2,3}, A_0 \subset A_1^3$
$\subset A_2^{2,3} \subset A_3^{2,3,4}, A_0 \subset A_1^3 \subset A_2^{3,4} \subset A_3^{1,3,4}, A_0 \subset A_1^3 \subset A_2^{3,4} \subset A_3^{2,3,4},$

$A_0 \subset A_1^4 \subset A_2^{1,4} \subset A_3^{1,2,4}, A_0 \subset A_1^4 \subset A_2^{1,4} \subset A_3^{1,3,4}, A_0 \subset A_1^4 \subset A_2^{2,4} \subset A_3^{1,2,4}, A_0 \subset A_1^4$
$\subset A_2^{2,4} \subset A_3^{2,3,4}, A_0 \subset A_1^4 \subset A_2^{3,4} \subset A_3^{1,3,4}, A_0 \subset A_1^4 \subset A_2^{3,4} \subset A_3^{2,3,4}.$

**Case II.** The chains in $FM_{2,3}(\mathcal{M})$ of the type

$$A_0 \subset A_1 \subset A_2 \subset A_4,$$

have the following 12 chains:

$A_0 \subset A_1^1 \subset A_2^{1,2} \subset A_4, A_0 \subset A_1^1 \subset A_2^{1,3} \subset A_4, A_0 \subset A_1^1 \subset A_2^{1,4} \subset A_4,$

$A_0 \subset A_1^2 \subset A_2^{1,2} \subset A_4, A_0 \subset A_1^2 \subset A_2^{2,3} \subset A_4, A_0 \subset A_1^2 \subset A_2^{2,4} \subset A_4,$

$A_0 \subset A_1^3 \subset A_2^{1,3} \subset A_4, A_0 \subset A_1^3 \subset A_2^{2,3} \subset A_4, A_0 \subset A_1^3 \subset A_2^{3,4} \subset A_4,$

$A_0 \subset A_1^4 \subset A_2^{1,4} \subset A_4, A_0 \subset A_1^4 \subset A_2^{2,4} \subset A_4, A_0 \subset A_1^4 \subset A_2^{3,4} \subset A_4.$

**Case III.** The chains in $FM_{2,3}(\mathcal{M})$ of the type

$$A_0 \subset A_1 \subset A_3 \subset A_4,$$

have the following 12 chains:

$A_0 \subset A_1^1 \subset A_3^{1,2,3} \subset A_4, A_0 \subset A_1^1 \subset A_3^{1,2,4} \subset A_4, A_0 \subset A_1^1 \subset A_3^{1,3,4} \subset A_4,$

$A_0 \subset A_1^2 \subset A_3^{1,2,3} \subset A_4, A_0 \subset A_1^2 \subset A_3^{1,2,4} \subset A_4, A_0 \subset A_1^2 \subset A_3^{2,3,4} \subset A_4,$

$$A_0 \subset A_1^3 \subset A_3^{1,2,3} \subset A_4, A_0 \subset A_1^3 \subset A_3^{1,3,4} \subset A_4, A_0 \subset A_1^3 \subset A_3^{2,3,4} \subset A_4,$$

$$A_0 \subset A_1^4 \subset A_3^{1,2,4} \subset A_4, A_0 \subset A_1^4 \subset A_3^{1,3,4} \subset A_4, A_0 \subset A_1^4 \subset A_3^{2,3,4} \subset A_4.$$

**Case IV.** The chains in $FM_{2,3}(\mathcal{M})$ of the type

$$A_0 \subset A_2 \subset A_3 \subset A_4,$$

have the following 12 chains:

$$A_0 \subset A_2^{1,2} \subset A_3^{1,2,3} \subset A_4, A_0 \subset A_2^{1,2} \subset A_3^{1,2,4} \subset A_4,$$

$$A_0 \subset A_2^{1,3} \subset A_3^{1,2,3} \subset A_4, A_0 \subset A_2^{1,3} \subset A_3^{1,3,4} \subset A_4,$$

$$A_0 \subset A_2^{1,4} \subset A_3^{1,2,4} \subset A_4, A_0 \subset A_2^{1,4} \subset A_3^{1,3,4} \subset A_4,$$

$$A_0 \subset A_2^{2,3} \subset A_3^{1,2,3} \subset A_4, A_0 \subset A_2^{2,3} \subset A_3^{2,3,4} \subset A_4,$$

$$A_0 \subset A_2^{2,4} \subset A_3^{1,2,4} \subset A_4, A_0 \subset A_2^{2,4} \subset A_3^{2,3,4} \subset A_4,$$

$$A_0 \subset A_2^{3,4} \subset A_3^{1,3,4} \subset A_4, A_0 \subset A_2^{3,4} \subset A_3^{2,3,4} \subset A_4.$$

**Case V.** The chains in $FM_{2,3}(\mathcal{M})$ of the type

$$A_1 \subset A_2 \subset A_3 \subset A_4,$$

have the following 24 chains:

$$A_1^1 \subset A_2^{1,2} \subset A_3^{1,2,3} \subset A_4, A_1^1 \subset A_2^{1,2} \subset A_3^{1,2,4} \subset A_4, A_1^1 \subset A_2^{1,3} \subset A_3^{1,2,3} \subset A_4, A_1^1 \subset A_2^{1,3} \subset A_3^{1,3,4} \subset A_4, A_1^1 \subset A_2^{1,4} \subset A_3^{1,2,4} \subset A_4, A_1^1 \subset A_2^{1,4} \subset A_3^{1,3,4} \subset A_4,$$

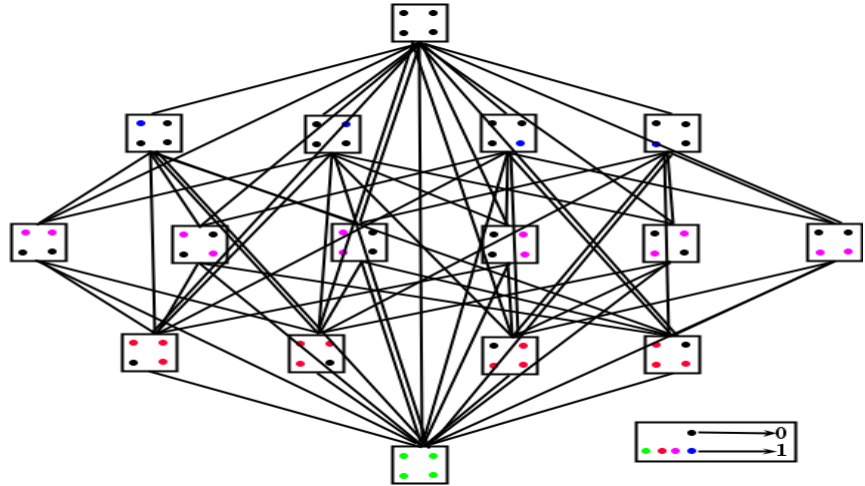

**Figure 2.** Graphical illustration of $FM_2(\mathcal{M})$.

$$A_1^2 \subset A_2^{1,2} \subset A_3^{1,2,3} \subset A_4, A_1^2 \subset A_2^{1,2} \subset A_3^{1,2,4} \subset A_4, A_1^2 \subset A_2^{2,3} \subset A_3^{1,2,3} \subset A_4, A_1^2$$
$$\subset A_2^{2,3} \subset A_3^{2,3,4} \subset A_4, A_1^2 \subset A_2^{2,4} \subset A_3^{1,2,4} \subset A_4, A_1^2 \subset A_2^{2,4} \subset A_3^{2,3,4} \subset A_4,$$

$$A_1^3 \subset A_2^{1,3} \subset A_3^{1,2,3} \subset A_4, A_1^3 \subset A_2^{1,3} \subset A_3^{1,3,4} \subset A_4, A_1^3 \subset A_2^{2,3} \subset A_3^{1,2,3} \subset A_4, A_1^3$$
$$\subset A_2^{2,3} \subset A_3^{2,3,4} \subset A_4, A_1^3 \subset A_2^{3,4} \subset A_3^{1,3,4} \subset A_4, A_1^3 \subset A_2^{3,4} \subset A_3^{2,3,4} \subset A_4,$$

$$A_1^4 \subset A_2^{1,4} \subset A_3^{1,2,4} \subset A_4, A_1^4 \subset A_2^{1,4} \subset A_3^{1,3,4} \subset A_4, A_1^4 \subset A_2^{2,4} \subset A_3^{1,2,4} \subset A_4, A_1^4$$
$$\subset A_2^{2,4} \subset A_3^{2,3,4} \subset A_4, A_1^4 \subset A_2^{3,4} \subset A_3^{1,3,4} \subset A_4, A_1^4 \subset A_2^{3,4} \subset A_3^{2,3,4} \subset A_4.$$

Thus, it is calculated from the above five cases that the number $f_{2,3}$ is

$$f_{2,3} = 24 + 12 + 12 + 12 + 24 = 84.$$

For all other particular cases, we will not provide a list of chains to avoid bulkiness. However, the graphical structure of $FM_2(\mathcal{M})$ have displayed in Figure 2. The observation through above examples and by analyzing Figure 1 and Figure 2, we come to conclude that the method of direct calculation doesn't work through manually (or becomes too complex) for other cases. Thus these problems motivate that we must find another method in order to determine the formula of $f_{n,k}$. Thus, we give the theorem as follows:

**Lemma 3.1.** *For any* $n \geq 1$ *and* $k \in \{0, 1, 2, \ldots, n^2\}$, *the number of* $FM_{n,k}(\mathcal{M})$'s *are given by the equality:*

$$f_{n,0} = \sum_{n_0=0}^{n^2} \binom{n^2}{n_0};$$

$$f_{n,1} = \sum_{n_0=0}^{n^2-1} \sum_{n_1=1}^{n^2-n_0} \binom{n^2}{n_0} \binom{n^2-n_0}{n_1};$$

$$f_{n,2} = \sum_{n_0=0}^{n^2-2} \sum_{n_1=1}^{n^2-n_0-1} \sum_{n_2=1}^{n^2-n_0-n_1} \binom{n^2}{n_0} \binom{n^2-n_0}{n_1} \binom{n^2-n_0-n_1}{n_2};$$

$$\vdots$$

$$f_{n,k} = \sum_{n_0=0}^{n^2-k} \sum_{n_1=1}^{n^2-n_0-k+1} \cdots \sum_{n_k=1}^{n^2-n_0-n_1-n_2-\cdots-n_{k-1}-k+k} \binom{n^2}{n_0} \binom{n^2-n_0}{n_1} \cdots \binom{n^2-n_0-n_1-n_2-\cdots-n_{k-1}}{n_k}.$$

*In particular, for* $k = n^2$ *we have the number* $f_{n,n^2}$ *of all distinct flags of crisp matrices of order* $n$ *is*

$$f_{n,n^2} = \prod_{k=1}^{n^2} \frac{(n^2-k+1)^{n^2-k+1}}{k!}.$$

**Proof.** We start the proof by making an auxiliary construction. It is known that all the crisp matrices of order $n$ are $A_0, A_1, A_2, \ldots, A_{n^2}$, where the number of entries in $A_k$, $|A_k| = \binom{n^2}{k}$ for $n \in \mathbb{N}$ and $k = 0, 1, 2, \ldots, n^2$. And these crisp matrices satisfy the following condition.

$$A_0 \subset A_1 \subset A_2 \subset \cdots \subset A_k, \text{ for } k = \overline{0, n^2}$$

Now
$FM_{n,0}(\mathcal{M}) = \{A_k \mid 0 \leq k \leq n^2\};$

$FM_{n,1}(\mathcal{M}) = \{A_i \subset A_k \mid 0 \leq i < k \leq n^2\};$

$\vdots$

$FM_{n,k}(\mathcal{M}) = \{A_{i_0} \subset A_{i_1} \subset A_{i_2} \subset \cdots \subset A_{i_k} \mid 0 \leq i_0 < i_1 < i_2 < \cdots < i_k \leq n^2\}.$

Next, our aim to find the cardinality of $FM_{n,k}(\mathcal{M})$ for $k = 0,1,2,\ldots,n^2$ and $\forall\, n \in \mathbb{N}$.

It can easily verify that $f_{n,0} = |FM_{n,0}(\mathcal{M})| = 2^{n^2}$. Let us consider $\Gamma$ be a chain of crisp matrices in $FM_{n,k}(\mathcal{M})$ for $1 \leq k \leq n^2$ as follows:

$$\Gamma : A_0 \subset A_1 \subset A_2 \subset \cdots \subset A_k.$$

Then one can obviously get the following $(k+1)$-dimensional vector as

$$\alpha = (|A_0|, |A_1|, |A_2|, \ldots, |A_k|).$$

For our convenience, let us call $\alpha$ the order vector of $\Gamma$.

Now consider

$$\Omega = \{\alpha \mid \alpha \text{ is an order vector of } \Gamma, \Gamma \in FM_{n,k}(\mathcal{M})\}$$

and for any $\alpha \in \Omega$,

$$FM_{n,k}^\alpha(\mathcal{M}) = \{\Gamma \in FM_{n,k}(\mathcal{M}) \mid \text{the order of } \Gamma \text{ is } \alpha\}.$$

Then it is very natural to see that

$$FM_{n,k}(\mathcal{M}) = \bigcup_{\alpha \in \Omega} FM_{n,k}^\alpha(\mathcal{M}),$$

and

$$FM_{n,k}^\alpha(\mathcal{M}) \cap FM_{n,k}^\beta(\mathcal{M}) = \emptyset \text{ if } \alpha \neq \beta.$$

Therefore,

$$f_{n,k} = \sum_{\alpha \in \Omega} |FM_{n,k}^\alpha(\mathcal{M})|$$

It is easy to notice that

$$\Omega = \left\{\alpha = \left(\binom{n^2}{i_0}, \binom{n^2}{i_1}, \binom{n^2}{i_2}, \ldots, \binom{n^2}{i_k}\right) \,\Big|\, 0 \leq i_0 < i_1 < i_2 < \cdots < i_k \leq n^2\right\}.$$

For the order vector $\alpha$ of $\Gamma \in FM_{n,k}(\mathcal{M})$ of length $k$, the number of choices for $1^{st}$ term, $2^{nd}$ term, $\cdots$, $k^{th}$ term of vector $\alpha$ are

$\binom{n^2}{n_0}, 0 \leq n_0 \leq n^2 - k$;

$\binom{n^2-n_0}{n_1}, 1 \leq n_1 \leq n^2 - n_0 - k + 1$;

$\vdots$

$\binom{n^2-n_0-n_1-n_2-\cdots n_{k-1}}{n_k}, 1 \leq n_k \leq n^2 - n_0 - n_1 - n_2 - \cdots n_{k-1}$.

Therefore, we can write exactly the formula $f_{n,k}$ as follows:

$f_{n,k} = \sum_{n_0=0}^{n^2-k} \sum_{n_1=1}^{n^2-n_0-k+1} \cdots \sum_{n_k=1}^{n^2-n_0-n_1-n_2-\cdots n_{k-1}-k+k} \binom{n^2}{n_0}\binom{n^2-n_0}{n_1}\cdots\binom{n^2-n_0-n_1-n_2-\cdots n_{k-1}}{n_k}$.

In particular, taking $k = n^2$ in $f_{n,k}$ we have

$$f_{n,n^2} = \prod_{k=0}^{n^2}\binom{n^2}{k} = \binom{n^2}{0}\binom{n^2}{1}\binom{n^2}{2}\cdots\binom{n^2}{n^2}$$

$$= \frac{(n^2)^{n^2}(n^2-1)^{n^2-1}(n^2-2)^{n^2-2}\cdots 1}{1!2!3!\cdots n^2!}$$

$$= \prod_{k=1}^{n^2} \frac{(n^2-k+1)^{n^2-k+1}}{k!}.$$

This completes the proof of the theorem. ∎

In order to understand and clarify the above Lemma 3.1, we will estimate the number $f_{2,2}$ in the following:

**Example 3.2.** Find the number $f_{2,2}$.

**Solution.** By using the explicit formula $f_{n,k}$ in Lemma 3.1. We shall find the number $f_{2,2}$ by taking $n = 2$ and $k = 2$. Then

$f_{2,2} = \sum_{n_0=0}^{2} \sum_{n_1=1}^{3-n_0} \sum_{n_2=1}^{4-n_0-n_1} \binom{4}{n_0}\binom{4-n_0}{n_1}\binom{4-n_0-n_1}{n_2} = \binom{4}{0}\binom{4}{1}\binom{3}{1} + \binom{4}{0}\binom{4}{1}\binom{3}{2} + \binom{4}{0}\binom{4}{1}\binom{3}{3} + \binom{4}{0}\binom{4}{2}\binom{2}{1} + \binom{4}{0}\binom{4}{2}\binom{2}{2} + \binom{4}{0}\binom{4}{3}\binom{1}{1} + \binom{4}{1}\binom{3}{1}\binom{2}{1} + \binom{4}{1}\binom{3}{1}\binom{2}{2} + \binom{4}{1}\binom{3}{2}\binom{1}{1} + \binom{4}{2}\binom{2}{1}\binom{1}{1} = 110$.

Hence the theorem is clarified.

In the following, the first two corollaries are obvious in the view of Lemma 3.1.

**Corollary 3.3.** *For each $n \in \mathbb{N}$ and $k = 0,1,\ldots,n^2$, the number of $FM_{n,k}^O(\mathcal{M})$'s are given by the equality:*

$f_{n,0}^O = 1$;

$$f_{n,1}^O = \sum_{n_1=1}^{n^2} \binom{n^2}{n_1};$$

$$f_{n,2}^O = \sum_{n_1=1}^{n^2-1} \sum_{n_2=1}^{n-n_1} \binom{n^2}{n_1}\binom{n^2-n_1}{n_2};$$

$\vdots$

$$f_{n,k}^O = \sum_{n_1=1}^{n^2-k+1} \sum_{n_2=1}^{n^2-n_1-k+2} \cdots \sum_{n_k=1}^{n^2-n_1-n_2-\cdots n_{k-1}} \binom{n^2}{n_1}\binom{n^2-n_1}{n_2}\cdots\binom{n^2-n_1-n_2-\cdots n_{k-1}}{n_k}.$$

**Corollary 3.4.** *For any non-zero positive integers $n$ and $0 \leq k \leq n^2$, the number of $FM_n^J(\mathcal{M})$'s, are given by the equality:*

$$f_{n,0}^J = 1;$$

$$f_{n,1}^J = \sum_{n_1=0}^{n^2-1} \binom{n^2}{n_1};$$

$$f_{n,2}^J = \sum_{n_1=0}^{n^2-2} \sum_{n_2=1}^{n^2-n_1-1} \binom{n^2}{n_1}\binom{n^2-n_1}{n_2};$$

$\vdots$

$$f_{n,k}^J = \sum_{n_1=0}^{n^2-k} \sum_{n_2=1}^{n^2-n_1-k+1} \cdots \sum_{n_k=1}^{n^2-n_1-n_2-\cdots n_{k-1}-1} \binom{n^2}{n_1}\binom{n^2-n_1}{n_2}\cdots\binom{n^2-n_1-n_2-\cdots n_{k-1}}{n_k}.$$

### 3.2. A Closed formula for $f_n$

The main objective of this subsection is to expose the closed summation formula $f_n$ for the number of all distinct fuzzy matrices of order $n$. By summing up all $f_{n,k}$'s we can easily obtain an explicit formula for $f_n$. So, the next result follows immediately from Lemma 3.1.

**Theorem 3.5.** *The number $f_n$ of all distinct fuzzy matrices of order $n$, $FM_n(\mathcal{M})$ is given by the following equality:*

$$f_n = \sum_{k=0}^{n^2} f_{n,k} =$$
$$\sum_{k=0}^{n^2}\left(\sum_{n_0=0}^{n^2-k}\sum_{n_1=1}^{n^2-n_0-k+1}\cdots\sum_{n_r=1}^{n^2-n_0-n_1-n_2-\cdots n_{r-1}-k+k}\binom{n^2}{n_0}\binom{n^2-n_0}{n_1}\cdots\binom{n^2-n_0-n_1-n_2-\cdots n_{k-1}}{n_k}\right),$$

*where $n \geq 1$ is arbitrary and fixed.*

Next, we are going to estimate the total number of terms in the expansion of the closed formula $f_n$ in the following proposition.

**Proposition 3.6.** *For any positive integers $n \geq 1$, the number of terms in the expansion of $f_n$ is $2^{n^2+1} - 1$.*

**Proof.** It is well-known that the number of all possible terms in the expression of $f_{n,k}$ is

$$\binom{n^2+1}{k+1}, k = 0,1,2,3, \ldots, n^2,$$

Therefore, the total number of terms in the expansion of the expression of the formula $f_n = \sum_{k=0}^{n^2} f_{n,k}$ is

$$\binom{n^2+1}{1} + \binom{n^2+1}{2} + \cdots + \binom{n^2+1}{n^2+1} = \sum_{k=1}^{n^2+1} \binom{n^2+1}{k},$$

by using Theorem 2.3, we can have

$$\sum_{k=1}^{n^2+1} \binom{n^2+1}{k} = 2^{n^2+1} - 1.$$

This completes the proof. ∎

We have the following two immediate straightforward consequence of Theorem 3.5.

**Corollary 3.7.** *For a fixed value $n \in \mathbb{Z}^+$ and $n \geq 1$, the number $f_n^O$ of all distinct O-rooted fuzzy matrices of order $n$ is given by the following equality:*

$$f_n^O = \sum_{k=0}^{n^2} f_{n,k}^O =$$
$$\sum_{k=0}^{n^2} (\sum_{n_1=1}^{n^2-k+1} \sum_{n_2=1}^{n^2-n_1-k+2} \cdots \sum_{n_k=1}^{n^2-n_1-n_2-\cdots n_{k-1}} \binom{n^2}{n_1}\binom{n^2-n_1}{n_2} \cdots \binom{n^2-n_1-n_2-\cdots n_{k-1}}{n_k}).$$

**Corollary 3.8.** *The number $f_n^J$ of all distinct J-rooted fuzzy matrices of order $n$, where $n \geq 1$, is given by the following equality:*

$$f_n^J = \sum_{k=0}^{n^2} f_{n,k}^J =$$
$$\sum_{k=0}^{n^2} (\sum_{n_1=0}^{n^2-k} \sum_{n_2=1}^{n^2-n_1-k+1} \cdots \sum_{n_k=1}^{n^2-n_1-n_2-\cdots n_{k-1}-1} \binom{n^2}{n_1}\binom{n^2-n_1}{n_2} \cdots \binom{n^2-n_1-n_2-\cdots n_{k-1}}{n_k}).$$

**Remark 3.9.** A fuzzy matrix of order 0, does not have a conventional meaning as there no elements that exist in the matrix. Thus, we simply take $f_0 = 1$ as an initial condition for the counting function $f_n$ for $n = 0$.

In the following, we construct the table for the number of $FM_{n,k}(\mathcal{M})$, $f_{n,k}$ (in Lemma 3.1) and the number of $FM_n(\mathcal{M})$, $f_n$ (in Theorem 3.5) for $0 \leq k \leq n^2, n \leq 3$.

**Table 1.** $f_n$ and $f_{n,r}$ for $0 \leq r \leq n^2, n \leq 3$

| n | $f_{n,0}$ | $f_{n,1}$ | $f_{n,2}$ | $f_{n,3}$ | $f_{n,4}$ | $f_{n,5}$ | $f_{n,6}$ | $f_{n,7}$ | $f_{n,8}$ | $f_{n,9}$ | $f_n$ |
|---|---|---|---|---|---|---|---|---|---|---|---|
| 0 | 1 | | | | | | | | | | 1 |
| 1 | 2 | 1 | | | | | | | | | 3 |
| 2 | 16 | 65 | 110 | 84 | 24 | | | | | | 299 |
| 3 | 512 | 19171 | 223290 | 1225230 | 3759840 | 6972840 | 8013600 | 5594400 | 2177280 | 362880 | 28349043 |

Thus, we have obtained the number of fuzzy matrices $f_n$ for $n \geq 0$, which, in turn, shall form a sequence $(f_n)_{n \geq 0}$ of natural numbers.

The initial first five terms for $n = 0,1,3,4,5$ of $(f_n)_{n \geq 0}$ are

$$1, 3, 299, 28349043, 21262618727925419.$$

**Table 2.** Number of fuzzy matrices $f_n$

| $n$ (A000027 in OEIS) | $n^2$ (A000290 in OEIS) | $2^{n^2}$ (A002416 in OEIS) | $f_n$ (Not available in OEIS) |
|---|---|---|---|
| 0 | 0 | 1 | 1 |
| 1 | 1 | 2 | 3 |
| 2 | 4 | 16 | 299 |
| 3 | 9 | 512 | 28349043 |
| 4 | 16 | 65536 | 21262618727925419 |
| 5 | 25 | 33554432 | . |
| 6 | 36 | 68719476736 | . |
| 7 | 49 | 562949953421312 | . |
| 8 | 64 | 18446744073709551616 | . |
| 9 | 81 | 2417851639229258349412352 | . |

## 4. Conclusions and Further Research

The investigate about the classification of fuzzy matrices is an important aspect of fuzzy matrix theory. In this paper, we have classified and counted all distinct fuzzy matrices of order $n$ by incorporating Murali's equivalence relation ([24], [33], [35]). For each $n \in \mathbb{N}$, listing the number of all distinct the fuzzy matrices of order $n$ will take a considerable amount of time due to enormous and the quick exponential growth of the number $f_n$. It will be a very complex task to calculate many more terms of the sequence $(f_n)$ using our counting technique. The classification of the counting problem can successfully be extended to some other special classes of fuzzy matrices. The study of fuzzy matrices can be effectively used in matrix theory, combinatorics as well as in geometry. It might be interesting to attack the classification of counting problem from different perspectives with more efficiently. This will surely indicate the way to further research.

There are two open problems according to this area of classification are as follows:

**Problem 4.1.** Is there any other explicit formula is to determine a formula for $f_n$ where $n$ is an arbitrary positive integer?

**Problem 4.2**. Write recurrence relation as well as generating function satisfied for the number $f_n$.